\documentclass{article}

\usepackage{amsmath} 
\usepackage{amsthm} 
\usepackage{graphics} 
\usepackage{hyperref} 
\usepackage{tikz}
\usepackage{amssymb}

\theoremstyle{plain}

\newtheorem{theorem}{{Theorem}}[section]

\newtheorem{definition}[theorem]{{Definition}}

\numberwithin{equation}{section} 

\author{Aaron Carl Smith}

\begin{document}


\title{Benford-Newcomb Subsequences for Fraud Detection}
\maketitle









\begin{abstract}
Benford's law is frequently used to evaluate the likihood that data is misrepresentative.  Typically statistical tests measure the likihood. Another method of employing Benford's law is to compare the frequency of leading digits to the probabilities of leading digits over a subset of the natural numbers.  This paper proposes using the probabilities of leading digits from uniform, natural numbers to establish interval criteria for when to look more closely into the possibility of misrepresentative data. 
\end{abstract}



\tableofcontents


\section{Introduction}Benford's law gives a probability distribution for the frequency of the leading-digit of natural numbers.  Simon Newcomb described the rule for decimal representation of natural numbers in 1881 \cite{MR1505286}, and Frank Benford generalized Newcomb's observations to any base in 1938 \cite{benford1938law}.  In 1995, Theodore Hill used the mantissa $\sigma$-algebra to further extend the leading-digit law to real numbers.  The mantissa $\sigma$-algebra consists of sets of numbers with the same coefficient in scientific notation after truncation \cite{MR1421567}.
\begin{definition}[Benford's Law] In base $b$, the probability that the leading digit of a real number is $k$ is given by
\begin{align}
P(k) = log_{b}(1 + \tfrac{1}{k}), \ k \in \{ 1,2,3,\dots, b-1\}.
\end{align}
\end{definition}
In decimal representation (base 10), the probabilities of each the leading digits are given by
\begin{align}
P(k) = log_{10}(1 + \tfrac{1}{k}), \ k \in \{1,2,3,4,5,6,7,8,9\},
\end{align} which approximately gives: \\ 
\begin{tabular}{|c|c|c|c|c|c|c|c|c|c|}
\hline
$k$ & $1$ & $2$ & $3$ & $4$ & $5$ & $6$ & $7$ & $8$ & $9$ \\
\hline
$P(k)$ & $0.301$ & $0.176$ & $0.125$ & $0.097$ & $0.079$ & $0.067$ & $0.058$ & $0.051$ & $0.046$ \\
\hline
\end{tabular}\\
The law goes further to say that the probability distribution of digits after the leading digit converges to uniform as the digit's position moves to the right \cite{benford1938law,MR1421567}.  Benford's law does not apply to several types of numeric data, such as identification numbers.
\section{Benford-Newcomb Subsequences}
Consider the map $f_b$ that sends natural numbers to their leading digits,
\begin{align}
f_b:\mathbb{N} \rightarrow \{1,2,3,\ldots,b-1 \}, x \mapsto floor(\tfrac{x}{b^{floor(log_b x)}}).
\end{align}        
Let $\mu_N$ be the uniform probability measure on $\mathbb{N}$ where $\mu_N(k) = \tfrac{1}{N} \ \forall \ k \in \{ 1,2,3,\ldots,N\}$.  Let's use $\mu_N$ to construct a probability measure of leading digits,
\begin{align}
P_{bN}(k) = \mu_N(\{x \in \mathbb{N} | f_b(x) = k \}).
\end{align}
For a fixed base $b$ and fixed leading digit $k$, consider the sequences $( P_{bN}(k)  )_{N = 1}^{\infty}$; in general these sequences do not converge.  The purpose of this paper is to propose using intervals of the form
\begin{align}
[\liminf_{N \rightarrow \infty}P_{bN}(k),\limsup_{N \rightarrow \infty}P_{bN}(k)]
\end{align} to identify possibly fraudulent data.  If a data set's frequency of leading digits, in base $b$ representation, is not contained in these intervals, then look further into the possibility of tamper data. For $N > b$, with respect to $N$ the local minimums are of the form
\begin{align}
P_{bN}(k) = \tfrac{1+b+b^2+\hdots+b^{\alpha - 1}}{kb^{\alpha} - 1}, \ N = kb^{\alpha} - 1
\end{align} and the local maximums are of the form
\begin{align}
P_{bN}(k) = \tfrac{1+b+b^2+\hdots+b^{\alpha}}{(k+1)b^{\alpha} - 1}, \ N = (k+1)b^{\alpha} - 1.
\end{align}
Thus if the frequencies of a data set's leading digits are not within
\begin{align}
[\tfrac{1}{k(b-1)},\tfrac{b}{(k+1)(b-1)}],
\end{align} further inquiry is called for.  The advantage of the interval method is that one may use it to quickly screen data.

\includegraphics[scale = 0.5]{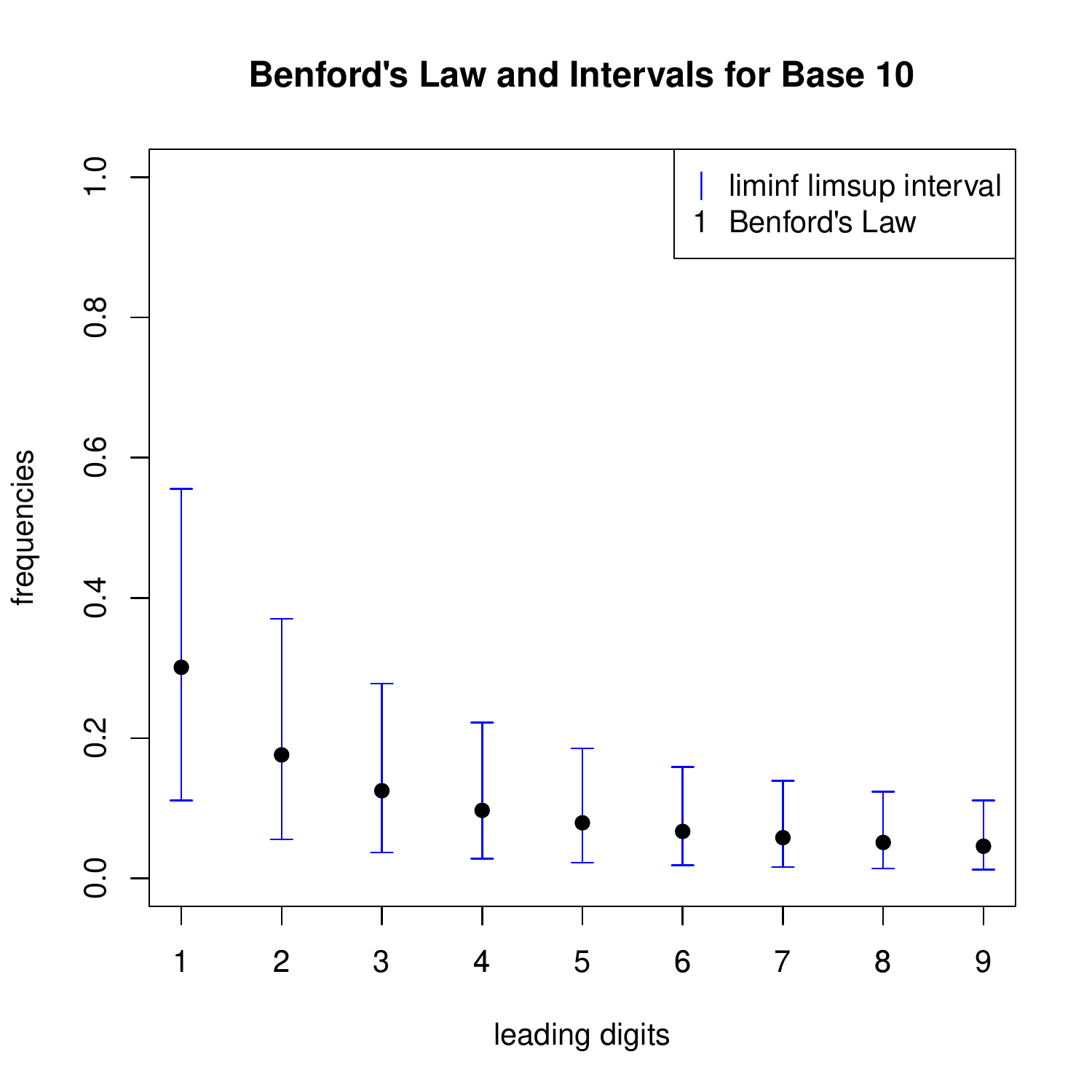}\includegraphics[scale = 0.45]{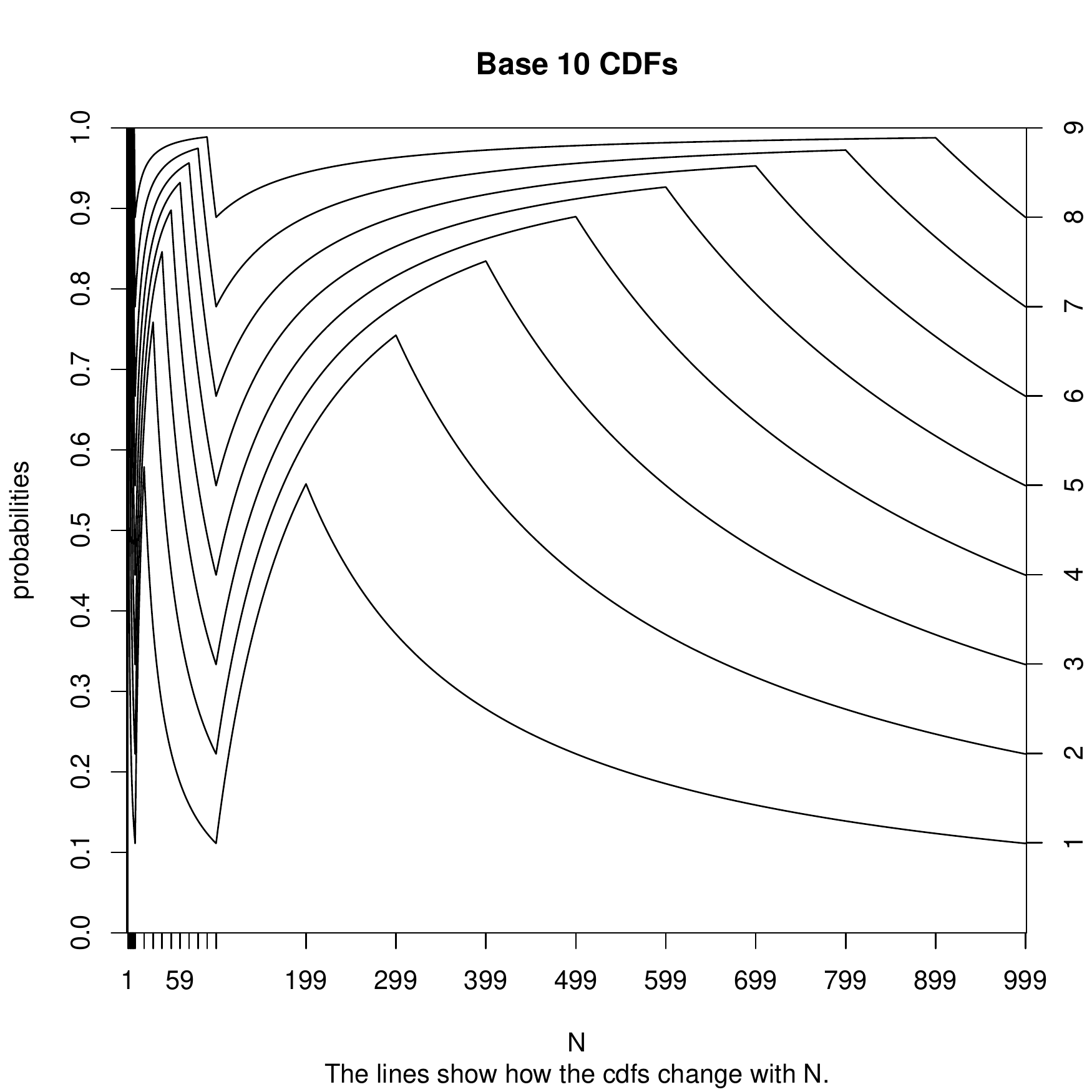} The figures were constructed with R \cite{RCoreTeam}.

\bibliographystyle{amsplain}
\bibliography{BenfordSubsequenceBibtex}

\end{document}